\documentclass{amsart}

\theoremstyle{definition}

\theoremstyle{remark}

\numberwithin{equation}{section}

\begin{document}

\title{An extension result of CR functions by a general Schwarz Reflection Principle}

\author{Chun Yin Hui}

\begin{abstract}
It is known that a real analytic CR function $f$ on a real
analytic, generic submanifold $M$ in $\mathbb{C}^N$ can be
holomorphically extended. A stronger result on a finite type, real analytic, generic submanifold
$M$ is found in which we
assume $f$ a continuous CR function with real analytic imaginary
part $Im(f)$. The idea is contained in a general Schwarz Reflection Principle
in one complex variable.\end{abstract}

\maketitle

\section{\textbf{A general Schwarz Reflection Principle}}
The classical Schwarz Reflection Principle can be stated as
follows:\\

\textbf{Theorem 1.1.} Suppose that $\Omega$ is a connected domain,
symmetric with respect to the real axis, and that
$L=\Omega\cap\mathbb{R}$ is an interval. Let
$\Omega^+=\{z\in\Omega:$ $Im(z)>0\}$. Suppose that $f\in
A(\Omega^+)$, a function holomorphic on $\Omega^+$ and that $Im(f)$
has a continuous extension to $\Omega^+ \cup L$ that vanishes on
$L$. Then there is a $F\in A(\Omega)$ such that $F=f$ in $\Omega^+$
and $F(z)=\overline{f(\overline{z})}$ in $\Omega-\Omega^+$.\\

Our generalized Schwarz Reflection Principle (Theorem 1.2) replaces the vanishing
of $Im(f)$ on $L$ by the real analyticity of $Im(f)$ on $L$.\\

\textbf{Theorem 1.2.} Let $D=D(0,1)$ be the open disc centered at
the origin with radius $1$, $D^+=\{z\in D$: $Im(z)>0\}$,
$L=D\cap\mathbb{R}$. Suppose that $f\in A(D^+)$ and $Im(f)\in
C(D^+\cup L)$ such that $v(x,0)=Im(f)|_L$ is real analytic at $0$
with radius of convergence $r$. Denote $D_r=D(0,r)$, then there
exists $F\in A(D^+\cup D_r)$ such that

\begin{equation}
F(z):= \left\{ \begin{array}{lll}
 f(z) &\mbox{for}& z\in D^+\\
 \overline{f(\overline{z})}+2iv(z,0) &\mbox{for}& z\in D_r-D_r^+
\end{array}\right.
\end{equation}

where $v(z,0)$ denotes a holomorphic function on $D_r$.\\

\textbf{Proof}: Since $v(x,0)=Im(f)|_L$ is real analytic at $0$ with
radius of convergence $r$, we can write $v(x,0)=\sum a_nx^n$,
$\forall |x|<r$. Thus we get a holomorphic function $v(z,0)=\sum
a_nz^n$ on $D_r$ by the complexification of the variable $x$. Now
$Im(f(z)-iv(z,0))|_{(-r,r)}\equiv0$, by the Schwarz Reflection
Principle, the function
\begin{equation}
g(z):= \left\{ \begin{array}{lll}
 f(z)-iv(z,0) &\mbox{for}& z\in D_r^+\\
 \overline{f(\overline{z})-iv(\overline{z},0)} &\mbox{for}& z\in D_r-D_r^+
\end{array}\right.
\end{equation}
is holomorphic on $D_r$. Thus the function
\begin{equation}
F(z):= g(z)+iv(z,0)= \left\{ \begin{array}{lll}
 f(z) &\mbox{for}& z\in D^+\\
 \overline{f(\overline{z})-iv(\overline{z},0)}+iv(z,0) &\mbox{for}& z\in D_r-D_r^+
\end{array}\right.
\end{equation}
is holomorphic on $D^+\cup D_r$ and is the desired extension. Since
$a_n$ are real, $F(z)=\overline{f(\overline{z})}+2iv(z,0)$, $\forall
z\in D_r-D_r^+$ which gives (1.1).
\hfill\rule{2.1mm}{2.1mm}\\

Theorem 1.2 leads to a reflection principle of harmonic functions.\\

\textbf{Corollary 1.3.} Let $D$, $D^+$, $L$ as above. Suppose that
$v(x,y)\in C(D^+\cup L)$ is harmonic in $D^+$ such that $v(x,0)$ is
real analytic at $0$ with radius of convergence $r$. Denote
$D_r=D(0,r)$, then there exists $V$ harmonic in $D^+\cup D_r$ such
that
\begin{equation}
V(x,y):= \left\{ \begin{array}{lll}
 v(x,y) &\mbox{for}& z\in D^+\\
 2Re(v(z,0))-v(x,-y) &\mbox{for}& z\in D_r-D_r^+
\end{array}\right.
\end{equation}
where $v(z,0)$ denotes a holomorphic function on $D_r$.\\

\textbf{Proof}: Since we can find $u(x,y)$ in $D^+$ such that
$u+iv\in A(D^+)$, we then apply Theorem 1.2 to get the extension
$V(x,y)$ in $D_r$,
$V(x,y)=Im(\overline{f(\overline{z})}+2iv(z,0))=2Re(v(z,0))-v(x,-y)$.
\hfill\rule{2.1mm}{2.1mm}\\

\textbf{Remark 1.4.} Since every real analytic curve in $\mathbb{C}$
is locally biholomorphic to the real line. The general Schwarz
Reflection Principle holds in the following case:\\

Let $S$ be a real analytic curve in $\mathbb{C}$ through $0$, $U$ a
small connected neighborhood of $0$ such that $U-S$ consists of two
components $U^+$, $U^-$. If $f\in A(U^+)$, the imaginary part
$Im(f)\in C(U-U^-)$ such that $Im(f)|_{U\cap S}$ is real analytic,
then $f$ can be extended to a function $F$ holomorphic in a
neighborhood of $0$ in $\mathbb{C}$.\\\\

\section{\textbf{Holomorphic extension of CR functions on a real analytic, generic CR submanifold in $\mathbb{C}^N$}}
We first introduce the necessary notations and definitions
needed in the sequel. We mainly follow [1]. For $Z\in\mathbb{C}^N$, we write $Z=(z_1,...,z_n)$ where $z_j=x_j+iy_j$; we write $\overline{Z}=(\overline{z_1},...,\overline{z_n})$, where $\overline{z_j}=x_j-iy_j$, the complex conjugate of $z_j$. We identify $\mathbb{C}^N$ with $\mathbb{R}^{2N}$, and denote a function $f$ on a subset of $\mathbb{C}^N$ as $f(x,y)$, or, by abuse of notation, as $f(Z,\overline{Z})$.\\

A smooth (real analytic)
real submanifold of $\mathbb{C}^N$ of codimension $d$ is a subset
$M$ of $\mathbb{C}^N$ such that $\forall p_0\in M$, there is a
neighborhood $U$ of $p_0$ and a smooth (real analytic) real
vector-valued function $\rho=(\rho_1,...,\rho_d)$ defined in $U$
such that
\begin{equation}
M\cap U=\{Z\in U: \rho(Z,\overline{Z})=0\},
\end{equation}
with differentials $d\rho_1,...,d\rho_d$ linearly independent in
$U$.\\

\textbf{Definition 2.1.} A real submanifold $M\subset \mathbb{C}^N$
is CR if the complex rank of the complex differentials
$\partial\rho_1,...,\partial\rho_d$ is constant for $p\in M$. It is
generic if $\partial\rho_1,...,\partial\rho_d$ are $\mathbb{C}$
linearly independent for $p\in M$.\\

Let
$M\subset\mathbb{C}^N$ be a germ of real analytic, generic, CR
submanifold of codimension $d$ at $p_0$, write $N=n+d$, $n\geq 1$.
After a local holomorphic change of coordinates, there exists $\Omega$, a sufficiently small open neighborhood of $0$
in $\mathbb{C}^{n+d}$, such that $M$ is given in $\Omega$ by
\begin{equation}
Im(w)=\phi(z, \overline{z}, Re(w)),
\end{equation}
with $z\in \mathbb{C}^n$, $w\in \mathbb{C}^d$, $\phi$ a
real-valued analytic function (power series) of $z, \overline{z}, Re(w)$
such that $\phi(0,0,0)=0$ and $d\phi(0,0,0)=0$. Such a choice of coordinates is called regular coordinates.\\

Now suppose $M$ is a real analytic, generic submanifold, locally
parametrized by the regular coordinates near $0\in M$. We
extend the parametrization
$\Psi(z,\overline{z},s)=(z,s+i\phi(z,\overline{z},s))\in\mathbb{C}^{n+d}$
to a local real analytic diffeomorphism
\begin{equation}
\widetilde{\Psi}(z,\overline{z},s+it)=(z,s+it+i\phi(z,\overline{z},s+it))\in\mathbb{C}^{n+d},
\end{equation}
defined in a neighborhood of $0\in\mathbb{R}^{2n+2d}$.\\

By the above notation, we have the following holomorphic extension
results, see [1].\\

\textbf{Proposition 2.2.} Let $M$ be a real analytic generic
submanifold of $\mathbb{C}^{n+d}$ of codimension $d$ in regular
coordinates near $0\in M$. If $h$ be a $CR$ function on $M$, then
$h$ extends holomorphically to a full neighborhood of $0$ if and
only if there exists $\epsilon>0$ such that for every
$z\in\mathbb{C}^n$, $|z|<\epsilon$, the function $s\longmapsto
h\circ\Psi(z,\overline{z},s)$ extends holomorphically to the open
set $\{s+it\in\mathbb{C}^d,$ $|s|<\epsilon$, $|t|<\epsilon\}$ in
such a way that the extension $H:=h\circ\Psi(z,\overline{z},s+it)$
is a bounded, measurable function of all its variables.\\

\textbf{Corollary 2.3.} Let $M\subset\mathbb{C}^N$ be a real
analytic generic submanifold and $f$ a CR function in a neighborhood
of $p\in M$. Then $f$ extends as a holomorphic function in a
neighborhood of $p$ in $\mathbb{C}^N$ if and only if $f$ is real
analytic in a neighborhood of $p$ in $M$.\\

According to Corollary 2.3, a CR fucntion $f=u+iv$ on a real
analytic, generic submanifold can be holomorphically extended to a
full neighborhood if and only if $u,v$ are real analytic on $M$. Our
main theorem, to some extent, reduces the case to $u$ continuous and
$v$ real analytic on $M$. We state it as follows:\\

(\textbf{Main theorem}) Let $M\subset\mathbb{C}^{n+d}$ be a real
analytic, generic submanifold. $0\in M$. In a neighborhood $\Omega$
of $0$ in $\mathbb{C}^{n+d}$, $M$ is given in regular coordinates
(2.2), $\rho:=Im(w)-\phi(z,\overline{z},Re(w))=0$,
$(z,w)\in\mathbb{C}^n\times\mathbb{C}^d$. Let $f=u+iv$ be a
holomorphic function defined in the wedge $W=W(\Omega,\rho,\Gamma)$,
with $\Gamma$ an open convex cone in $\mathbb{R}^d$. If $f$ extends
continuously up to the edge $M\cap\Omega$ with $v=Im(f)$ is real
analytic on the edge $M\cap\Omega$. Then there exists a holomorphic
function $F$ in a neighborhood $\widetilde{\Omega}$ of $0$ in
$\mathbb{C}^{n+d}$
such that $F|_{W\cap \widetilde{\Omega}}=f$.\\

Let $M$ be a generic submanifold
of $\mathbb{C}^{n+d}$ of codimension $d$ and $p_0\in M$. Let
$\rho=(\rho_1,...,\rho_d)$ be a defining functions of $M$ near $p_0$
and $\Omega$ a small neighborhood of $p_0$ in $\mathbb{C}^{n+d}$ in
which $\rho$ is defined. If $\Gamma$ is an open convex cone with
vertex at the origin in $\mathbb{R}^d$, we define
\begin{equation}
W(\Omega,\rho,\Gamma):=\{Z\in\Omega: \rho(Z,\overline{Z})\in\Gamma\}.
\end{equation}
The above set is an open subset of $\mathbb{C}^{n+d}$ whose boundary
contains $M\cap\Omega$. Such a set is called a wedge of edge $M$ in the direction $\Gamma$ centered at $p_0$.\\

The following shows that $W(\Omega,\rho,\Gamma)$ is in a sense
independent of the choice of $\rho$, which will allow us to change
defining functions freely. See [1].\\

\textbf{Proposition 2.4.} Let $\rho$ and $\rho'$ be two defining
functions for $M$ near $p_0$, where $M$ and $p_0$ are as above. Then
there is a $d\times d$ real invertible matrix $B$ such that for
every $\Omega$ and $\Gamma$ as above the following holds. For any
open convex cone $\Gamma_1\subset\mathbb{R}^d$ with $\{Bx:$
$x\in\Gamma_1\}\cap S^{d-1}$ relatively compact in $\Gamma\cap
S^{d-1}$ (where $S^{d-1}$ denotes the unit sphere in
$\mathbb{R}^d$), there exists $\Omega_1$, an open neighborhood of
$p_0$ in $\mathbb{C}^{n+d}$, such that
\begin{equation}
W(\Omega_1,\rho',\Gamma_1)\subset W(\Omega,\rho,\Gamma).
\end{equation}

Below is an Edge-of-the-Wedge Theorem, see [2], [3]. Since the
construction of the holomorphic extension in this theorem is
essential to the proof of the main theorem, so we supply a proof of
the Edge-of-the-Wedge
Theorem from [2], p.157-159.\\

\textbf{Theorem 2.5.} (\textbf{Edge-of-the-Wedge Theorem}) If
$\Gamma\subset\mathbb{R}^d$ is an open convex cone, $R>0$, $d\geq
2$. Let $V=\Gamma\cap B(0,R)$. Let $E\subset\mathbb{R}^d$ be a
nonempty neighborhood of $0$. Define $W^+\subset\mathbb{C}^d$,
$W^-\subset\mathbb{C}^d$ by
\begin{equation}
W^+=E+iV, ~~W^-=E-iV
\end{equation}
Then there exists a fixed neighborhood $U$ of $0\in\mathbb{C}^d$
such that the following property holds: For any continuous function
$g:W^+\cup W^-\cup E\longrightarrow\mathbb{C}$ that is holomorphic
on $W^+\cup W^-$, there is a holomorphic $G$ on $U$ such that
$G|_{U\cap(W^+\cup W^-\cup E)}=g$.\\

\textbf{Proof}: After composition with a linear isomorphism $A$ in
$\mathbb{C}^d$ with real coefficients, we may assume that
$\{(it_1,...,it_d)\in\mathbb{C}^d:$ $t_j>0,$ $j=1,...,d\}\subset
A^{-1}(i\Gamma)$, there exists $B(0,R')\subset\mathbb{R}^d$ with
$R'>6\sqrt{d}$ and $E'=\{s\in\mathbb{R}^d:$ $|s_j|<6,$ $j=1,...,d\}$
such that $E'+iB(0,R')\subset A^{-1}(E+iB(0,R))$. So we reduce the case to the following:\\

\noindent
$E=\{s\in\mathbb{R}^d:$ $|s_j|<6,$ $j=1,...,d\}$,\\
$V=\{t\in\mathbb{R}^d:$ $0<t_j<6,$ $j=1,...,d\}$,\\
$W^+=E+iV$,~~$W^-=E-iV$.\\

Let $U=D^d(0,1)$. If $g:W^+\cup W^-\cup E\longrightarrow\mathbb{C}$
is continuous and $g$ is holomorphic on $W^+\cup W^-$, then there is
a holomorphic $G$ on $U$ such that $G|_{U\cap(W^+\cup W^-\cup
E)}=g$. Let $c=\sqrt{2}-1$ and define
\begin{equation}
\varphi:\overline{D}^2(0,1)\longrightarrow\mathbb{C}
~~~~~~~~~~(w,\lambda)\longmapsto\frac{w+\lambda/c}{1+c\lambda w}
\end{equation}
Then
\begin{equation}
Im\varphi(w,\lambda)=\frac{(1-|\lambda|^2)Im(cw)+(1-|cw|^2)Im(\lambda)}{c|1+c\lambda w|^2}
\end{equation}
Notice that
\begin{enumerate}
\item[(a)] $sgn(Im\varphi)=sgn(Im\lambda)$ if $|\lambda|=1$ or $w\in\mathbb{R}$;
\item[(b)] $\varphi(w,0)=w$;
\item[(c)] $|\varphi(w,\lambda)|\leq(1+1/c)(1-c)<6$;
\item[(d)] By (a), the function\\
$\Phi:D^d(0,1)\times\overline{D}\longrightarrow\mathbb{C}^d$ with
$\Phi(w,\lambda)=(\varphi(w_1,\lambda),...,\varphi(w_d,\lambda))$\\
satisfies $\Phi(w,e^{i\theta})\in$ Dom$(g)$ for all
$0\leq\theta<2\pi$, $w\in D^d(0,1)$.
\end{enumerate}

So by (d), we may define
\begin{equation}
G(w)=\frac{1}{2\pi}\int_0^{2\pi}g(\Phi(w,e^{i\theta}))d\theta, ~~~w\in D^d(0,1).
\end{equation}
We claim that this is what we seek. First, $G$ is holomorphic by an
application of Morera's theorem. Next, for fixed $s\in E\cap
D^d(0,1)$, the function $g(\Phi(s,\cdot))$ is continuous on
$\overline{D}$ and holomorphic on $D-\mathbb{R}$ by (a),(c). Again,
by Morera, $g(\Phi(s,\cdot))$ is holomorphic on all of $D$. It
follows by Mean Value Property and (b) that
\begin{equation}
G(s)=\frac{1}{2\pi}\int_0^{2\pi}g(\Phi(s,e^{i\theta}))d\theta=g(\Phi(s,0))=g(s);
\end{equation}
hence $G=g$ on $E\cap D^d(0,1)$. If $s+it\in\mathbb{C}^d$ is fixed, $|s+it|<1/2$, $t>0$, then the function
\begin{equation}
\xi\longmapsto G(s+\xi t)-g(s+\xi t)
\end{equation}
is holomorphic for $\xi\in\mathbb{C}$ small and $\equiv 0$ when $\xi$ is real.
If follows that $G\equiv g$ on $W^+\cup W^-\cup E$.\\

If we return to our original case, the desired extension of $g$ is
\begin{equation}
G(w)=\frac{1}{2\pi}\int_0^{2\pi}g\circ A(\Phi(A^{-1}(w),e^{i\theta}))d\theta, w\in A^{-1}(D^d(0,1)).
\end{equation}
since it agrees with our $g$ in some open sets. $U$ is taken to be
$A^{-1}(D^d(0,1))$
\hfill\rule{2.1mm}{2.1mm}\\

\textbf{Theorem 2.6.} (\textbf{Main theorem}) Let
$M\subset\mathbb{C}^{n+d}$ be a real analytic, generic submanifold.
$0\in M$. In a neighborhood $\Omega$ of $0$ in $\mathbb{C}^{n+d}$,
$M$ is given by regular coordinates,
$\rho:=Im(w)-\phi(z,\overline{z},Re(w))=0$,
$(z,w)\in\mathbb{C}^n\times\mathbb{C}^d$. Let $f=u+iv$ be a
holomorphic function defined in the wedge $W=W(\Omega,\rho,\Gamma)$,
where $W(\Omega,\rho,\Gamma)$ is as above defined with $\Gamma$ an
open convex cone in $\mathbb{R}^d$. If $f$ extends continuously up
to the edge $M\cap\Omega$ with $v=Im(f)$ is real analytic on the
edge $M\cap\Omega$. Then there exists a holomorphic function $F$ in
a neighborhood $\widetilde{\Omega}$ of $0$ in $\mathbb{C}^{n+d}$
such that $F|_{W\cap \widetilde{\Omega}}=f$.\\

Before proving the theorem, we need a lemma which ensures the
existence of a wedge with edge as an open subet of
$\mathbb{R}^{2n+d}$ contained in
$\widetilde{\Psi}^{-1}(W)$. $\widetilde{\Psi}$ is the real analytic diffeomorphism introduced in (2.3).\\

\textbf{Lemma 2.7.} There exists a neighborhood $U_1$ of $0$ in
$\mathbb{R}^{2n+2d}$, a neighborhood $Q_1$ of $0$ in
$\mathbb{R}^{2n+d}$, an open convex cone $\Gamma'$ with vertex at
$0$ in $\mathbb{R}^d$ such that $\widetilde{\Psi}$ is a real
analytic diffeomorphism from $U_1$ to
$\widetilde{\Psi}(U_1)\subset\Omega$ and
$\widetilde{\Psi}((Q_1\times\Gamma')\cap U_1)\subset W$.\\

\textbf{Proof of the lemma}: Take a neighborhood $U_1$ of $0$ in
$\mathbb{R}^{2n+2d}$ such that $\widetilde{\Psi}$ is the real
analytic diffeomorphism from $U_1$ to
$\widetilde{\Psi}(U_1)\subset\Omega$.
$\widetilde{\Psi}(z,s+it)=(z,s+it+i\phi(z,\overline{z},s+it))=(\widetilde{z},\widetilde{w})$,
we can write
$t=\Theta(\widetilde{z},\overline{\widetilde{z}},\widetilde{w},\overline{\widetilde{w}})$
where $\Theta$ is $\mathbb{R}^d$ valued, real analytic in all its
variables. Note that
$t:=\Theta(\widetilde{z},\overline{\widetilde{z}},\widetilde{w},\overline{\widetilde{w}})$
can be taken as a set of defining functions for $M$ in
$\widetilde{\Psi}(U_1)$ because $\widetilde{\Psi}$ is a real
analytic diffeomorphism. By Proposition 2.4, there exists a
neighborhood of $0$ in $\mathbb{C}^{n+d}$(we may take it as a subset
of $\Psi(U_1)$ and still denote it as $\Psi(U_1)$), a convex open
cone $\Gamma'$ such that
\begin{equation}
W(\Psi(U_1),t,\Gamma')\subset W=W(\Omega,\rho,\Gamma),
\end{equation}
Take $Q_1=U_1\cap\mathbb{R}^{2n+d}$.
Let $(\widetilde{z},\widetilde{w})=\widetilde{\Psi}(x,t)$, where
\begin{equation}
(x,t)=(x_1,...,x_{2n+d},t_1,...,t_d)\in(Q_1\times\Gamma')\cap U_1.
\end{equation}
Since $(t_1,...,t_d)\in\Gamma'$, thus
$\Theta(\widetilde{z},\overline{\widetilde{z}},\widetilde{w},\overline{\widetilde{w}})\in\Gamma'$,
$(\widetilde{z},\widetilde{w})\in W(\Psi(U_1),t,\Gamma')\subset
W(\Omega,\rho,\Gamma)$. So we have
$\widetilde{\Psi}((Q_1\times\Gamma')\cap U_1)\subset W$.

\hfill\rule{2.1mm}{2.1mm}\\

\textbf{Proof of the theorem}: $f|_{M\cap\Omega}$ is a continuous CR
function on $M\cap\Omega$. By the lemma, the function
$f\circ\widetilde{\Psi}(z,\overline{z},s+it)=f(z,\overline{z},s+it+i\phi(z,\overline{z},s+it))$
is defined on $\{(z,s+it)\in(Q_1\times\Gamma')\cap U_1\}$, note that
it is holomorphic in the variables $w=s+it$. To use Proposition 2.2,
it suffices to show that $f\circ\widetilde{\Psi}$ can be
continuously extended to a full neighborhood of $0$ in
$\mathbb{R}^{2n+2d}$ such that the extension is holomorphic in
variables
$w=s+it$.\\

Now, denote the restriction of the function
$v\circ\widetilde{\Psi}(z,s+it)$ on $Q_1$ to be
$v(z,\overline{z},s)$. Note that the real analytic diffeomorphism
(holomorphic in $s+it$)$\widetilde{\Psi}$ can be viewed as the
flattening of the totally real submanifolds
$\widetilde{\Psi}(z,\overline{z},s)$ for fixed $z$. By the real
analyticity of $v$ and $\Psi$, we can treat $v(z,\overline{z},s)$ as
a power series expansion at $0\in\mathbb{R}^{2n+d}$ and complexify
the variable $s$ to $w$. Denote the new function to be
$v(z,\overline{z},w)$ defined in a product neighborhood
$U_2=U_3\times U_4\subset\mathbb{C}^n\times\mathbb{C}^d$ of $0$.
Note that $v(z,\overline{z},w)$ is continuous in $U_2$ and
holomorphic in variables $w$.\\

Following the idea of Theorem 1.2, the function
$g(z,\overline{z},w)$
\begin{equation}
:= \left\{ \begin{array}{lll}
 f\circ\widetilde{\Psi}(z,\overline{z},w)-iv(z,\overline{z},w) &\mbox{for}& (z,w)\in(Q_1\times(\Gamma'\cup\{0\}))\cap
 U_2\\
 \overline{f\circ\widetilde{\Psi}(z,\overline{z},\overline{w})-iv(z,\overline{z},\overline{w}}) &\mbox{for}& (z,w)\in(Q_1\times-\Gamma')\cap U_2
\end{array}\right.
\end{equation}
is holomorphic in variables $w=s+it$ for fixed $z$ and continuously up to $Q_1\cap U_2$.\\

We shall apply the construction of holomorphic extension (2.12) in the
proof of Theorem 2.5 to $g(z,\overline{z},s+it)$ for every fixed $z$
to obtain a continuous extension in all variables. According to the
proof, we have a linear isomorphism $A$ in $\mathbb{C}^d$ with real
coefficients and
$\Phi:D^d(0,1)\times\overline{D}\longrightarrow\mathbb{C}^d$. Now we
extend $A$ to a linear isomorphism in $\mathbb{C}^{n+d}$ which maps
$(z,w)\in\mathbb{C}^{n+d}$ to $(z,A(w))$, still denote it as $A$. We
also extend $\Phi$ to a mapping which takes $(z,w,\lambda)\in
U_3\times D^d(0,1)\times\overline{D}$ to
$(z,\Phi(w,\lambda))\in\mathbb{C}^n\times\mathbb{C}^d$, still denote
it as $\Phi$. Since the choice of $A$ depends only on $U_4$ and the
cone $\Gamma'$, so we have the extension
\begin{equation}
G(z,\overline{z},w)=\frac{1}{2\pi}\int_0^{2\pi}g\circ A(\Phi(A^{-1}(z,w),e^{i\theta}))d\theta,
~\forall(z,w)\in A^{-1}(U_3\times D^d(0,1))
\end{equation}
by Theorem 2.5 and Lemma 2.7.\\

From the construction (2.12), the extension $G(z,\overline{z},s+it)$
is continuous on $A^{-1}(U_3\times D^d(0,1))=U_3\times
A^{-1}(D^d(0,1))=U_3\times U_5$ and holomorphic in variables $w$.
Thus, the function
$F(z,\overline{z},s+it):=G+iv(z,\overline{z},s+it)$ is a continuous
extension of $f\circ\widetilde{\Psi}(z,\overline{z},s+it)$ in
$U_3\times U_5$, holomorphic in $w=s+it$. Thus by Proposition 2.2,
$f$ can be holomorphically extended to a full neighborhood
$\widetilde{\Omega}$ of $0$ in $\mathbb{C}^{n+d}$.
\hfill\rule{2.1mm}{2.1mm}\\

The following theorem states that minimality is a sufficient
condition for holomorphic extension of all CR functions from a
generic submanifold $M$ in $\mathbb{C}^N$ into an open wedge. See
[1], [4].\\

\textbf{Theorem 2.8.} Let $M$ be a generic submanifold of
$\mathbb{C}^{n+d}$ of codimension $d$ and $p_0\in M$. If $M$ is
minimal at $p_0$, then for every open neighborhood $U$ of $p_0$ in
$M$ there exists a wedge $W$ with edge $M$ centered at $p_0$ such
that every continuous CR function in $U$ extends holomorphically to
the wedge $W$.\\

Since a real analytic CR submanifold $M$ in $\mathbb{C}^N$ is finite
type at $p_0\in M$ if and only if it is minimal at $p_0$. Thus,
together with Theorem 2.6, we have Corollary 2.9.\\

\textbf{Corollary 2.9.} Let $M\subset\mathbb{C}^N$ be a real
analytic, generic CR submanifold, finite type at $p_0\in M$. If
$f=u+iv$ is a continuous CR function defined in a neighborhood of
$p_0$ in $M$ with $v$ real analytic. Then $f$ can be holomorphically
extended to a full neighborhood of $p_0$ in $\mathbb{C}^N$. (Or $u$
is also real analytic in a neighborhood of $p_0$ in $M$ by Corollary
2.3)\\

\textbf{Corollary 2.10.} Let $M$ be a connected, real analytic,
generic CR submanifold in $\mathbb{C}^N$. Assume $M$ of finite type at its every point. If
$f,g$ are continuous CR functions on $M$ such that $Im(f)=Im(g)$ on
$M$ and $f(p_0)=g(p_0)$ at $p_0\in M$. Then $f=g$ on $M$.\\

\textbf{Proof}: $f-g$ is a continuous CR function on $M$.
$Im(f-g)\equiv 0$ is real analytic. By finite type condition at every point and Corollary 2.9., $f-g$ can be extended
to a holomorphic function $H$ defined in a connected neighborhood $U$ of $M$ in $\mathbb{C}^N$.\\

We shall show that $H$ is constant on $U$. Suppose not, without loss of generality, we assume that $\frac{\partial H}{\partial z_N}$ is not constant zero on $U$.
Define subsets in $U$ as follows:
\begin{equation}
Q:=\{Z\in U: \frac{\partial H}{\partial z_N}(Z)=0\}
\end{equation}
\begin{equation}
I:=\{Z\in U: Im(H)(Z)=0\}
\end{equation}
Clearly, $M\subset I$. Since $M$ is a generic submanifold and $Q$ is an analytic variety in $U$, $M\not\subset Q$.
Let $p_1\in M-Q$, then $\frac{\partial H}{\partial z_N}(p_1)\neq0$. Thus we take a small neighborhood $U_1$ of $p_1$ in $U$ such that $I\cap U_1$ is a smooth hypersurface and $U_1\ni(z_1,...,z_{N-1},z_N)\mapsto(z_1,...,z_{N-1},H(Z))$ is a biholomorphism. This implies $I\cap U_1$ is biholomorphic equivalent to the real hyperplane which is of infinite type.  However, $p_1\in M\subset I$, $M$ is of finite type at $p_1$, hence $I\cap U_1$ is of finite type at $p_1$. This contradiction shows that $H$ is constant on $U$. By assumption, $H(p_0)=0$ and $U$ is connected. Hence $H\equiv 0$ on $U$. So $f\equiv g$ on $M$. \hfill\rule{2.1mm}{2.1mm}\\

\textbf{Definition 2.11.} A CR submanifold of the form
\begin{equation}
M=\{(z,s+it)\in\mathbb{C}^n\times\mathbb{C}^d; t=\phi(z,\overline{z})\}
\end{equation}
where $\phi:\mathbb{C}^n\longmapsto\mathbb{R}^d$ is smooth with
$\phi(0)=0$ and $d\phi(0)=0$ is called rigid.\\

\textbf{Corollary 2.12.} Let $M$ be a real analytic, generic
submanifold, $0\in M$, given by $\{(z,s+it)\in
U\subset\mathbb{C}^n\times\mathbb{C}^d:$
$t=\phi(z,\overline{z},s)\}$. Let $M'$ be a real analytic, generic,
rigid submanifold, $0\in M'$, given by $\{(z',s'+it')\in
U'\subset\mathbb{C}^{n'}\times\mathbb{C}^{d'}:$
$t'=\phi'(z',\overline{z'})\}$. If
$H=(f_1,...,f_{n'},g_1,...,g_{d'})$ is a holomorphic mapping defined
in the wedge $W=W(U,\rho,\Gamma)$ with $\Gamma$ an open convex cone.
If $H$ extends continuously to the edge $M\cap U$ with $H(0)=0$,
$H(M)\subset M'$ and $(f_1,...,f_{n'})$ can be holomorphically
extended in a neighborhood of $0$, then $H$ can be holomorphically
extended in a neighborhood of
$0$.\\

\textbf{Proof}: It suffices to show $Im(g_i)$ is real analytic for
all $i$. By the assumption, $f_1,...,f_{n'}$ can be holomorphically
extended in a neighborhood of $0$. By Corollary 2.3, they are real
analytic on $M$, hence their real and imaginary parts are real
analytic on $M$. Since $\phi'_i$ is real analytic in a neighborhood
of $0$ in $\mathbb{C}^{n'}$ and $Im(g_i)=\phi'_i(Re(f),Im(f))$,
hence $Im(g_i)$ are real analytic on $M$ for all $i$. Thus by
Theorem 2.6, $g_i$ can be holomorphically extended to a neighborhood
of $0$, we get the desired result.
\hfill\rule{2.1mm}{2.1mm}\\

\bibliographystyle{amsplain}

\bigskip
\bigskip

\noindent Einstein Institute of Mathematics,
Hebrew University, Givat Ram,
Jerusalem 91904, Israel\\
Email: \texttt{pslnfq@gmail.com, chhui@umail.iu.edu}
\end{document}